\documentclass[11pt,leqno,english]{article}

\usepackage{ShVa644}

\begin{document}

%\InputText{}{TitleAndAbstract}
%begin(TitleAndAbstract)===========================================
%\RequireCommands{}{system}

\title{\protect
	On inverse $\gamma$-systems and
	the number of \Lan [\lambda]-equivalent,
	non-isomorphic models for $\lambda$ singular
}
%EndTitle

\author{
	Saharon Shelah
	\thanks{Thanks to GIF for its support of this research and
	also to University of Helsinki for funding a visit of the
	first author to Helsinki in August 1996. Pub. No. 644.} \\
	\and
	Pauli V\"{a}is\"{a}nen
	\thanks{The second author wishes to thank Tapani Hyttinen
	under whose supervision he did his share of the paper.}
}
%EndAuthors

\date{\today}

\maketitle

\begin{abstract}%
Suppose $\lambda$ is a singular cardinal of uncountable cofinality
$\kappa$. For a model \M of cardinality $\lambda$, let $\No(\M)$
denote the number of isomorphism types of models \N of cardinality
$\lambda$ which are \Lan [\lambda]-equivalent to \M. In \cite {Sh189}
Shelah considered inverse $\kappa$-systems \System [A] of abelian
groups and their certain kind of quotient limits \GrF [A]. In
particular Shelah proved in \cite [Fact 3.10] {Sh189} that for every
cardinal $\mu$ there exists an inverse $\kappa$-system \System [A]
such that \System [A] consists of abelian groups having cardinality at
most $\mu^\kappa$ and $\CardGrF [A] = \mu$. Later in \cite [Theorem
3.3] {Sh228} Shelah showed a strict connection between inverse
$\kappa$-systems and possible values of \No \Note{under the assumption
that $\theta^\kappa < \lambda$ for every $\theta < \lambda$}: if
\System [A] is an inverse $\kappa$-system of abelian groups having
cardinality $< \lambda$, then there is a model \M such that $\Card \M
= \lambda$ and $\No (\M) = \CardGrF [A]$. The following was an
immediate consequence \Note{when $\theta^\kappa < \lambda$ for every
$\theta < \lambda$}: for every nonzero $\mu < \lambda$ or $\mu =
\lambda^\kappa$ there is a model $\M_\mu$ of cardinality $\lambda$
with $\No(\M_\mu) = \mu$.  In this paper we show: for every nonzero
$\mu \leq \lambda^\kappa$ there is an inverse $\kappa$-system \System
[A] of abelian groups having cardinality $< \lambda$ such that
$\CardGrF [A] = \mu$ \Note {under the assumptions $2^\kappa < \lambda$
and $\theta^{<\kappa} < \lambda$ for all $\theta < \lambda$ when $\mu
> \lambda$}, with the obvious new consequence concerning the possible
value of \No. Specifically, the case $\No(\M) = \lambda$ is possible
when $\theta^\kappa < \lambda$ for every $\theta < \lambda$.%
\footnote{%
}%
 \end{abstract}

%end(TitleAndAbstract)=============================================

\begin{SECTION} {-} {Introduction} {Introduction}

%\SInput
%begin(Section_Introduction)=======================================
%\RequireCommands{}{system}

Suppose $\lambda$ is a cardinal. For a model \M we let \Card \M denote
the cardinality of the universe of \M. When \M and \N are models of
the same vocabulary and they satisfy the same sentences of the
infinitary language \Lan [\lambda], we write $\M \LEquiv [\lambda]
\N$. For any model \M of cardinality $\lambda$ we define $\No (\M)$
to be the cardinality of the set
 \[
	\Set [\big] {\Quotient \N \Isomorphic} {
		\Card \N = \lambda \And
		\N \LEquiv [\lambda] \M 
	}, %EndSet
 \]
 where \Quotient \N \cong is the equivalence class of \N under the
isomorphism relation. Our principal purpose is to study the possible
values of $\No (\M)$ for models \M of singular cardinality with
uncountable cofinality.

% ------- History -------------------------------------------------

When \M is countable, $\No(\M) = 1$ by \cite {Scott}. This result
extends to structures of cardinality $\lambda$ when $\lambda$ is a
singular cardinal of countable cofinality \cite {Chang}.

If $V = L$, $\lambda$ is an uncountable regular cardinal which is not
weakly compact, and \M is a model of cardinality $\lambda$, then $\No
(\M)$ has either the value 1 or $2^\lambda$. For $\lambda = \aleph_1$
this result was first proved in \cite {Palyutin}. Later in \cite
{Sh129} Shelah extended this result to all other regular non-weakly
compact cardinals. The possibility $\No (\M) = \aleph_0$ is consistent
with $\ZFC + \GCH$ in case $\lambda =\aleph _{1}$, as remarked in
\cite {Sh129}. The values $\No (\M) \in \omega \Minus \Braces {0,1}$
are proved to be consistent with $\ZFC + \GCH$ in the forthcoming
paper of the authors \cite {ShVa646} \Note {number 646 in Shelah's
publications}.

The case \M has cardinality of a weakly compact cardinal is dealt with
in \cite {Sh133} by Shelah. The result is that for $\kappa$ weakly
compact there is for every $1 \leq \mu \leq \kappa$ a model $\M_\mu$
such that $\No(\M_\mu) = \mu$. There is in preparation by the authors
a paper where the question for $\kappa$ weakly compact is revisited.

% -------- This paper -----------------------------------------------

The case \M is of singular cardinality $\lambda$ with uncountable
cofinality $\kappa$ was first treated in \cite {Sh189}, where the
relations of \M have infinitely many places. Later in \cite {Sh228}
Shelah improved the result by showing that if $\theta^\kappa <
\lambda$ for every $\theta < \lambda$ and $0 < \mu < \lambda$ then
$\No(\M) = \mu$ is possible for a model \M having cardinality
$\lambda$ and relations of finitely many places only. The main idea in
those papers was to transform the problem of possible values of $\No
(\M)$ into a question concerning possible cardinalities of ``quotient
limit'' \GrF [A] of an inverse system \System [A] of groups \cite
[Theorem 3.3] {Sh228}:

\begin{THEOREM}{ModelsAndSystems}%
\Note {$\lambda$ cardinal with $\lambda > \Cf \lambda = \kappa >
\aleph_0$}%
If $\theta^\kappa < \lambda$ for every $\theta < \lambda$ and \System
[A] is an inverse $\kappa$-system of abelian groups having cardinality
$< \lambda$, then there is a model \M of cardinality $\lambda$ \Note
{with relations having finitely many places only} such that $\No (\M)
= \CardGrF [A]$.
\end{THEOREM}

Actually the groups in \cite [Theorem 3.3] {Sh228} are not limited to
be abelian. However, abelian groups suffice for the present purposes.

The recent paper fills a gap left open since the paper \cite
{Sh228}. We present a uniform way to construct inverse $\kappa$-system
of abelian groups having a quotient limit of desired cardinality. The
most important new case is that the cardinality of a quotient limit
can be $\lambda$ for some inverse system \Note {in other cases, where
the result below can be applied, the Singular Cardinal Hypothesis
fails}. The result of this paper is:

\begin{THEOREM}{System}%
\Note {$\lambda$ cardinal with $\lambda > \Cf \lambda = \kappa >
\aleph_0$}%
For every nonzero $\mu \leq \lambda$ there is an inverse
$\kappa$-system $\System [A] = \SystemSeq [\kappa]$ of abelian groups
satisfying that $\Card {G_i} < \lambda$ for every $i < \kappa$ and
$\CardGrF [A] = \mu$. The same conclusion holds also for the values
$\lambda < \mu \leq \lambda^\kappa$ under the assumption that
$2^\kappa < \lambda$ and $\theta^{<\kappa} < \lambda$ for every
$\theta < \lambda$.
\end{THEOREM}

So the general method used here to find new possibilities for the
values of $\No(\M)$ is the same as in \cite {Sh228}. As an immediate
consequence of the last theorem we get:

\begin{THEOREM}{Models}%
Suppose $\lambda$ is a singular cardinal of uncountable cofinality
$\kappa$. For each nonzero $\mu \leq \lambda^\kappa$ there is a model
\M \Note {with relations having finitely many places only} satisfying
$\Card \M = \lambda$ and $\No (\M) = \mu$, provided that
$\theta^\kappa < \lambda$ for every $\theta < \lambda$.
\end{THEOREM}

We give all necessary definitions concerning inverse $\kappa$-systems
\System [A] of abelian groups and their special kind of quotient
limits \GrF [A] in the next section.

%end(Section_Introduction)=========================================

\end{SECTION}

\begin{SECTION} {-} {Preliminaries} {Preliminaries}

%\SInput
%begin(Section_Preliminaries)=======================================
%\RequireCommands{}{system}
 
\begin{DEFINITION}{General_system}
 Suppose $\gamma$ is a limit ordinal and for every $i < j < \gamma$,
\G i is a group and \Hom [i,j] is a homomorphism from \G j into \G
i. The family $\System [A] = \SystemSeq [\gamma]$ is called an inverse
$\gamma$-system when the equation $\Comp {\Hom [i,j]} {\Hom [j,k]} =
{\Hom [i,k]}$ holds for every $i < j < k < \gamma$. As in \cite
{Sh189} we assume that all the groups $G_i$, $i < \gamma$, are
additive abelian groups.

To simplify our notation we make an agreement that the letters $i$,
$j$, $k$, and $l$ always denote ordinals smaller than $\gamma$. Hence
``for all $i < j$'' means ``for all ordinals $i$ and $j$ with $i < j <
\gamma$'' and so on.

The main objects of our study are the following two sets:
 \begin{align}
	\Gr [A] &= \Set [\Big] { \Seq {\GS [i,j] a} {i < j < \gamma} }
	{\GS [i,j] a \in \G i \Text{ and for all} k > j,\ 
	\GS [i,k] a = \GS [i,j] a + \Hom [i,j] (\GS [j,k] a)}; \notag\\
	\Fact [A] &= \Set [\Big] { \Seq {\GS [i,j] a} {i < j < \gamma} }
	{\ForSome \YS \in \Prod [k < \gamma] {\G k},\ 
	\GS [i,j] a = \YS [i] - \Hom [i,j] (\YS [j])}. \notag
 \end{align}
 We consider \Gr [A] and \Fact [A] as additive abelian groups where
the group operation $+$ and the unit element \Zero are pointwise
defined. The factor group \GrF [A] is well-defined since $\Fact [A]
\Subset \Gr [A]$ by the requirements $\Comp {\Hom [i,j]} {\Hom [j,k]}
= {\Hom [i,k]}$ for all $i < j < k$. For any inverse $\gamma$-system
\System [A], the group \GrF [A] is called the quotient limit of
\System [A].
 \end{DEFINITION}

\begin{DEFINITION}{I-sets}
 We let \Ind [\gamma] be the set \Set {(i,j) \in \gamma \times \gamma}
{i < j}. For every subset $I$ of \Ind we define
 \[
	\First I = \Set {i < \gamma} {(i,j) \in I \ForSome j < \gamma}
 \]
 and for each $i \in \First I$,
 \[
	\Pro I i = \Set {j < \gamma} {(i,j) \in I}.
 \]
 We also say that
 \begin{itemize}

\item
 $I$ is cobounded if $\gamma \Minus \First I$ and $\gamma \Minus \Pro
I i$, for all $i \in \First I$, are bounded subsets of $\gamma$;

\item
 $I$ is coherent if \First I is unbounded in $\gamma$ and for every $i
\in \First I$, $\Pro I i = \First I \Minus (i+1)$;

\item
 $I$ is eventually coherent if it is unbounded and for every $i \in
\First I$, $\First I \Minus \Pro I i$ is a bounded subset of $\gamma$.

\end{itemize}
\end{DEFINITION}

\Remark Suppose $I$ is an eventually coherent subset of \Ind [\gamma]
and $S$ is a subset of \First I. If $\Card S < \Cf \gamma$, then
$\First I \Minus (\BigInter [i \in S] {\Pro I i})$ is a bounded subset
of $\gamma$. If $S$ is unbounded in $\gamma$, then $I \Inter (S \times
S)$ is an eventually coherent subset of $I$.

In \cite [Claim 1.12] {Sh228} Shelah proved \Note {note the remark
given after the following lemma} that if two sequences \GS a and \GS b
from \Gr [A] agree on a coherent set of indices, then \FactEquivGS [A]
a b. The following slight improvement of this condition has an
essential role in the proof of \Theorem {System}.

\begin{LEMMA}{Eventually_coherent}
 Suppose \System [A] is an inverse $\gamma$-system, and $\GS a, \GS b
\in \Gr [A]$. Then \FactEquivGS [A] a b holds if there is an
eventually coherent subset $I$ of \Ind such that $\GS [i,j] a = \GS
[i,j] b$ for all $(i,j) \in I$.

\begin{Proof}
 We shall need an eventually coherent subset $J$ of $I$ having the
property that \Seq {\Pro J i} {i \in \First J} is a decreasing chain
of end segments of \First J. Let $S$ be an unbounded subset of $I$
having the order type \Cf \gamma. Define a subset $J$ of $I$ by
$\First J = S$ and for all $j \in S$,
 \[
	\Pro J j = S \Inter
	\BigInter [i \in S \Inter (j+1)]
		{\Par{\Pro I i \Minus (i^*+1)}},
 \]
 where $i^*$ is the supremum of the bounded subset $\First I \Minus
\Pro I i$ of $\gamma$. The set $J$ is well-defined since $I$ is
eventually coherent and $\Card [\big] {S \Inter (j+1)} < \Cf \gamma$
for all $j < \gamma$. Now $J$ is also eventually coherent, and
furthermore, for all $i \in \First J$, $\Pro J i = S \Minus \Min {\Pro
J i}$ and for all $j \in \First J \Minus i$, $\Min {\Pro J i} \leq
\Min {\Pro J j}$.

Define for every $i < \gamma$, $i'$ to be \Min [\big] {\First J \Minus
(i+1)} and $i'' = \Min {\Pro J {i'}}$. Then the following are
satisfied for all $i < j$:
 \begin{itemize}

\item $i< i' < i''$, $j < j' < j''$, $i' \leq j'$, $i'' \leq j''$, and also
$i' < j''$;

\item $j'' \in \First I$ and $(i',i''), (j', j''), (i',j'') \in I$.

\end{itemize}

Since \GS a and \GS b are in \Gr we have
\ARRAY[lll]{
	\GS [i,j''] a = \GS [i,j] a + \Hom [i,j] (\GS [j,j''] a), \\
	\GS [i,j''] b = \GS [i,j] b + \Hom [i,j] (\GS [j,j''] b).
}%EndArray
 Therefore the following equations hold:
 \begin{align}
 \AEQUATION{i,j}{
	\GS [i,j] a - \GS [i,j] b &=
		(\GS [i,j''] a - \GS [i,j''] b)
		- \Par {\Hom [i,j] (\GS [j,j''] a)
		- \Hom [i,j] (\GS [j,j''] b)}
 }\\
 \AEQUATION{}{
	&= (\GS [i,j''] a - \GS [i,j''] b)
		- \Hom [i,j] (\GS [j,j''] a - \GS [j,j''] b).
 }
\end{align}
 Because of $i < i' < j''$ we also have that
 \ARRAY[lll]{
	\GS [i,j''] a = \GS [i,i'] a + \Hom [i,i'] (\GS [i',j''] a), \\
	\GS [i,j''] b = \GS [i,i'] b + \Hom [i,i'] (\GS [i',j''] b).
 }
 Since $(i',j'') \in I$, $\GS [i',j''] a = \GS [i',j''] b$
holds. Hence we get
 \EQUATION{i,j''}{%
	\GS [i,j''] a - \GS [i,j''] b = \GS [i,i'] a - \GS [i,i'] b.%
 }%
 Moreover, $i < i' < i''$ yields
 \ARRAY[lll]{
	\GS [i,i''] a = \GS [i,i'] a + \Hom [i,i'] (\GS [i',i''] a), \\
	\GS [i,i''] b = \GS [i,i'] b + \Hom [i,i'] (\GS [i',i''] b).
 }
 Now $(i',i'') \in I$ implies that $\GS [i',i''] a = \GS [i',i''] b$,
and consequently
 \[
	\GS [i,i'] a - \GS [i,i'] b = \GS [i,i''] a - \GS [i,i''] b.
 \]
 This equation together with \Equation {i,j} and \Equation {i,j''}
implies that for all $i < j$
 \[
	\GS [i,j] a - \GS [i,j] b =
	(\GS [i,i''] a - \GS [i,i''] b)
		- \Hom [i,j] (\GS [j,j''] a - \GS [j,j''] b).
 \]
 So the sequence $\YS = \Seq {\GS [i,i''] a - \GS [i,i''] b} {i <
\gamma} \in \Prod [i < \gamma] {G_i}$ exemplifies that $\GS a - \GS b
\in \Fact [A]$, and we have \FactEquivGS [A] a b.
 \end{Proof}
\end{LEMMA}

\Remark In \cite [Claim 1.12] {Sh228} the groups of an inverse system
\System [A] need not to be abelian groups. Hence instead of the factor
group \GrF [A] a partition \Quotient {\Gr [A]} {\approx_{\System [A]}}
with a special kind of equivalence relation $\approx_{\System [A]}$
were considered there. However, it is straightforward to prove, by
means of the preceding proof, also the more general case of \Lemma
{Eventually_coherent} where ``equivalent modulo \Fact [A]'' is
replaced by $\approx_{\System [A]}$.

In the next section we shall need a notion of a tree, so we shortly
describe our notation.

\begin{DEFINITION}{Tree}
 Suppose $T = \Structure {T, \Tleq}$ is a tree of height $\gamma$.
For every $i < \gamma$, \Level i is the \Th i level of the tree. When
$i < j < \gamma$ and $\eta \in \Level j$, then \Res \eta i denotes the
unique element $\nu \in \Level i$ for which $\nu \Tleq \eta$
holds. For each $i < \gamma$ and $\nu \in \Level i$, \ProLevel j \nu
is the set \Set {\eta \in \Level j} {\nu \Tleq \eta}. The set of all
$\gamma$-branches of $T$, i.e., the set \Set {t \in \Prod [i < \gamma]
{\Level i}} {\ForAll i < j,\ t(i) \Tleq t(j)}, is denoted by \Br
[\gamma].
 \end{DEFINITION}

%end(Section_Preliminaries)=========================================

\end{SECTION}

\begin{SECTION} {-} {System} {\protect%
	The inverse $\gamma$-system of free $R$-modules}

%\SInput
%begin(Section_System)=======================================
%\RequireCommands{}{system}

In this section we define special kind of inverse $\gamma$-systems
\System and prove a result concerning cardinalities of their quotient
limit \GrF \Note {\Conclusion {System}}. A direct consequence of the
result will be \Theorem {System}.

\begin{DEFINITION}{System}
 Suppose $\gamma$ is a limit ordinal, $R$ is a ring, and $T$ is a tree
of height $\gamma$. We define an inverse $\gamma$-system $\System =
\SystemSeq$ by the following stipulations:
  \begin{ITEMS}

\ITEM {Modules}
 for each $i < \gamma$, $G_i$ is the $R$-module freely generated by
\Set {\Gen [\nu,l]} {\nu \in \Level i \And i < l < \gamma};

\ITEM {Hom}
 for every $i < j < \gamma$, \Hom [i,j] is the homomorphism from $G_j$
into $G_i$ determined by the values $\Hom [i,j] (\Gen [\eta,l]) =
\ProGen [i,l] \eta - \ProGen [i,j] \eta$, for all $\eta \in \Level j$
and \GreInd l j. \Note {It is easy to check that the equations $\Hom
[i,k] = \Comp {\Hom [i,j]} {\Hom [j,k]}$ are satisfied for all $i < j
< k$}.

 \end{ITEMS}
 We consider \Gr, \Fact, and \GrF as $R$-modules where the operations
$+, \Times$, and the unit element \Zero for addition are pointwise
defined.

 For each $t \in \Br$, we define \GenS t to be the sequence \Seq
{\GenT [i,j] t} {i < j < \gamma}. Directly by the definitions of $G_i$
and \Hom [i,j], \GenS t belongs to \Gr for every $t \in \Br$. We let
\BrGen be the submodule of \Gr generated by the elements \GenS t, $t
\in \Br$. When \Br is empty \BrGen is the trivial submodule \Braces
\Zero.
 \end{DEFINITION}

\Remark Each $G_i$ is nonempty when $T$ has
height $\gamma$. Hence \Prod [i < \gamma] {G_i} is nonempty, and also
 \[
	\Fact = \Set [\Big] {
		\Seq {\YS [i] - \Hom [i,j] (\YS [j])}
			{i < j < \gamma}
	}{\YS \in \Prod [i < \gamma] {G_i}}
 \]
 is nonempty. So $\Gr \Superset \Fact$ is nonempty for every ring $R$
and tree $T$ of height $\gamma$.

Observe also that the inverse $\gamma$-system \System is the same as
used in \cite [Claim 3.8] {Sh189} when $R$ is the trivial ring \Braces
{0,1} and $T$ consists of $\mu$ many disjoint $\gamma$-branches. So
the proof given in this section offers an alternative proof for \cite
[Claim 3.8] {Sh189}, and even more information, namely that \CardGrF
must be exactly $\mu$ not only $\geq \mu$.

\begin{DEFINITION}{Sum}
 Suppose $\GS a \in \Gr$ and $i < j < \gamma$. By the definition of \G
i and the requirement $\GS [i,j] a \in \G i$, we define $\Coff [i,j] a
{\nu,l}$ for $\nu \in \Level i$ and $l > i$, to be the coefficients
from $R$ \Note {with only finitely many of them nonzero} which satisfy
the equation
 \[
	\GS [i,j] a = \GSSum {i,j} a {\nu,l} \CoffGen.
 \]
 The finite set \Set {(\nu,l) \in \Level i \times \Par {\gamma \Minus
(i+1)}} {\Coff [i,j] a {\nu,l} \not= 0} is called the support of \GS
[i,j] a, and it is denoted by \Supp [i,j] a.

Suppose $S$ is a subset of $\gamma$, $e \in G_i$, and $e_{\nu,l} \in
R$ for every $\nu \in \Level i$ and $l > i$ are elements such that
 \[
	e = \RSum {\nu \in \Level i} {\GreInd l i}
		{e_{\nu,l} \Times \Gen [\nu,l]}.
 \]
 Then we write \Res e S for the following element of $G_i$:
 \[
	\RSum {\nu \in \Level i} {l \in S\Minus(i+1)}
		{e_{\nu,l} \Times \Gen [\nu,l]}.
 \]
 \end{DEFINITION}

The following simple lemma has an important corollary.

\begin{LEMMA}{Simple}
 \begin{ITEMS}

\ITEM{Hom}
 The restriction \Res {\Hom [i,j] (e)} j equals 0 for every $i<j$
and $e \in G_j$.

\ITEM{Fact}
 For every $\GS a \in \Fact \Minus \Braces \Zero$, there are $i < j <
\gamma$ such that $\GSRes [i,j] a j \not= 0$.

\end{ITEMS}

\begin{Proof}
 \ProofOfItem {Hom}
 Straightforwardly by the definitions of $G_j$ and \Hom [i,j].

 \ProofOfItem {Fact}
 By the definition of \Fact, let $\YS \in \Prod [i < \gamma] {G_i}$ be
such that for all $i<j$, $\GS [i,j] a = \YS [i] - \Hom [i,j] (\YS
[j])$. In addition to that let $y^i_{\nu,l} \in R$, for $i < \gamma$,
$\nu \in \Level i$ and \GreInd l i, be such that
 \[
	\YS [i] = \RSum {\nu \in \Level i} {\GreInd l i}
	{y^i_{\nu,l} \Times \Gen [\nu,l]}.
 \]
 Since $\GS a \not= \Zero$ there must be $i < \gamma$ with $\YS [i]
\not= 0$. Define $j$ to be $\min \Set {\GreInd l i} {y^i_{\nu,l} \not=
0 \ForSome \nu \in \Level i} +1$. Then $\Res {\YS [i]} j$ is nonzero
and because $\Res {\Hom [i,j] (\YS [j])} j = 0$, we have $\GSRes [i,j]
a j = \Res {\YS [i]} j - \Res {\Hom [i,j] (\YS [j])} j = \Res {\YS
[i]} j \not= 0$.
 \end{Proof}
\end{LEMMA}

\begin{COROLLARY}{Independent}
 The elements \GenS t, $t \in \Br$, are independent over \Fact, i.e.,
\[
	\BrGen \Inter \Fact = \Braces \Zero.
\]
 Hence \System satisfies $\CardGrF \geq \Card [\big] \BrGen$.

\begin{Proof}
Directly by the definition of \GenS t, $\GenS [i,j] t = \GenT [i,j] t$
and hence $\GenSRes [i,j] t j = 0$, for all $t \in \Br$ and $i<j$.  So
for any nonzero $\GS a = \Sum [1 \leq m \leq n] {d_m \Times \GenS
{t_m}}$, where $n < \omega$, $d_m \in R \Minus \Braces 0$, and $t_m
\in \Br$, the restrictions $\GSRes [i,j] a j$ are equal to 0 for all
$i<j$. So by the preceding lemma \GS a can not be in \Fact.
\end{Proof}
\end{COROLLARY}

Next we derive equations of weighty significance.

\begin{LEMMA}{Eq}
 Suppose $\GS b \in \Gr$ and $i < j < k < \gamma$. Then the following
equations are satisfied for all $\nu \in \Level i$:
 \begin{alignat}{2}
 \AEQUATION{<j}{
	\Coff [i,k] b {\nu,l} &= \Coff [i,j] b {\nu,l}
	&\qquad &\When\ \BetInd i l j;
 }\\
 \AEQUATION{=j}{
	\Coff [i,k] b {\nu,j} &=
	\Coff [i,j] b {\nu,j} - \GSSumResT {j,k} b {\eta,l} \Coff \nu;
	& &
 }\\
 \AEQUATION{>j}{
	\Coff [i,k] b {\nu,l} &=
	\Coff [i,j] b {\nu,l} +
	\GSSumResT [] {j,k} b {\eta,l} \Coff \nu
	&\qquad &\When\ \GreInd l j.
 }
\end{alignat}

\begin{Proof}
 By dividing the sum into groups we get that
 \begin{equation*}\begin{split}
	\GS [i,j] b &= \GSSum {i,j} b {\nu,l} \CoffGen \\
	&= \Sum [\nu \in \Level i] {\Par[\Big]{%
		  \Sum [\BetInd i l j] {\CoffGen [i,j] b {\nu,l}}
		+ \CoffGen [i,j] b {\nu,j}
		+ \Sum [\GreInd l j]    {\CoffGen [i,j] b {\nu,l}}
	}}.%EndParSum
 \end{split}\end{equation*}
 Similarly the following equation is satisfied,
\[
	\GS [i,k] b
	= \Sum [\nu \in \Level i] {\Par[\Big]{%
		  \Sum [\BetInd i l j] {\CoffGen [i,k] b {\nu,l}}
		+ \CoffGen [i,k] b {\nu,j}
		+ \Sum [\GreInd l j]    {\CoffGen [i,k] b {\nu,l}}
	}}.%EndParSum
\]
 From the definition of \Hom [i,j] we may infer that
 \ARRAY[lll]{
	\Hom [i,j] (\GS [j,k] b)
	&=& \RSum {\eta \in \Level j} {\GreInd l j}
		{\Coff[j,k] b {\eta,l} \Times
		\Hom [i,j] (\Gen [\eta,l])} \\
	&=& \RSum {\eta \in \Level j} {\GreInd l j}
		{\Coff[j,k] b {\eta,l} \Times
		(\ProGen [i,l] \eta - \ProGen [i,j] \eta)} \\
	&=& \RSum {\eta \in \Level j} {\GreInd l j}
		{\Coff [j,k] b {\eta,l} \Times \ProGen [i,l] \eta}
	  - \RSum {\eta \in \Level j} {\GreInd l j}
		{\Coff [j,k] b {\eta,l} \Times \ProGen [i,j] \eta} \\
	&=& \Sum [\nu \in \Level i] {\Par[\Big]{%
		\Sum [\GreInd l j] {%
		    \Par[\big]{%
			\Sum [\eta \in \ProLevel j \nu]
			{\Coff [j,k] b {\eta,l}}
		    } \Times \Gen [\nu,l]
		}%EndSum
		-   \Par[\big] {%
			\GSSumResT {j,k} b {\eta,l} \Coff \nu
		    }
		    \Times \Gen[\nu,j]
	 }%End\Par[\bigg]
	}.%EndSum
}%EndArray
 So the equations \Equation {<j}, \Equation {=j}, and \Equation {>j}
for all $i<j<k$ follow by comparing the coefficients of each generator
\Gen [\nu,l] in the equation $\GS [i,k] b = \GS [i,j] b + \Hom[i,j]
(\GS [j,k] b)$.
 \end{Proof}
\end{LEMMA}

\begin{LEMMA} {Res}
 Suppose $\GS a \in \Gr$.
 \begin{ITEMS}

\ITEM {Res}
 For all $i<j<k$, \GSRes [i,j] a j = \GSRes [i,k] a j.

\ITEM {Finite_union} \Note{$\Cf \gamma > \aleph_0$}\ 
 For every $i < \gamma$, the union \BigUnion [i < j < \gamma]
{\ResSupp [i,j] a j} is of finite cardinality \Note {where $\ResSupp
[i,j] a j = \Supp [i,j] a \Inter (\Level i \times j)$ of course}.

\ITEM {Cobounded} \Note{$\Cf \gamma > \aleph_0$}\ 
 There is $\GS b \in \Gr$ satisfying the following conditions:
 \begin{itemize}

\item \FactEquivGS a b,

\item I = \Set {(i,j) \in \Ind [\gamma]} {\GSRes [i,j] b j = 0} is
cobounded \Note {in fact $\First I = \gamma$}, and

\item for every $(i,j) \in I$, $\GS [i,j] b = \GSRes [i,j] b {\Braces
j}$.

 \end{itemize}
 \end{ITEMS}

\begin{Proof}
 \ProofOfItem {Res}
 The claim holds directly by \EquationOfLemma {Eq} {<j}.

\ProofOfItem {Finite_union}
 Suppose the union is infinite. Since $\Cf \gamma > \aleph_0$ there is
some $k < \gamma$ for which already \BigUnion [j < k] {\ResSupp [i,j]
a j} is infinite. By \Item {Res}, $\ResSupp [i,j] a j \Subset \Supp
[i,k] a$ for each $j < k$. Consequently $\BigUnion [j < k] {\ResSupp
[i,j] a j} \Subset \Supp [i,k] a$ contrary to the finiteness of \Supp
[i,k] a.

\ProofOfItem {Cobounded}
 By \Item {Res} and \Item {Finite_union} there must be for every $i <
\gamma$ a bound $i^* \in \gamma \Minus (i+1)$ such that for every $j
\geq i^*$, $\GSRes [i,i^*] a {i^*} = \GSRes [i,j] a j$. Define an
element $\GS c \in \Fact$ by
 \[
	\GS [i,j] c = \GSRes [i,i^*] a {i^*}
		- \Hom [i,j] (\GSRes [j,j^*] a {j^*}),
 \]
 for all $i < j$. Let \GS b be $\GS a - \GS c$. Then \FactEquivGS a b
and for every $i < \gamma$ and $j \geq i^*$
 \ARRAY[lll]{
	\GS [i,j] b &=& \GS [i,j] a - \GS [i,j] c \\
	&=& \GSRes [i,j] a {(\gamma\Minus j)} + \GSRes [i,j] a j
		- \GSRes [i,i^*] a {i^*}
		+ \Hom [i,j] (\GSRes [j,j^*] a {j^*}) \\
	&=& \GSRes [i,j] a {(\gamma\Minus j)}
		+ \Hom [i,j] (\GSRes [j,j^*] a {j^*}).
 }%EndArray
 It follows from \ItemOfLemma {Simple} {Hom} that $\GSRes [i,j] b j =
0$ for all $i < \gamma$ and $j \geq i^*$, and thus $I$ is cobounded.

Now suppose, contrary to the last claim in \Item {Cobounded}, that
$\Coff [i,j] b {\nu,l} \not= 0$ for some $i < \gamma$, $j \geq i^*$,
$\nu \in \Level i$, and \GreInd l j. Let $k$ be $\Max {i^*, j^*,
l+1}$. Then both \GSRes [i,k] b k and \GSRes [j,k] b k are 0. By
\EquationOfLemma {Eq} {>j} the following equation holds:
 \[
	\GSSumResT [] {j,k} b {\eta,l} \Coff \nu
	= \Coff [i,k] b {\nu,l} - \Coff [i,j] b {\nu,l}.
 \]
 Since $\Coff [i,j] b {\nu,l} \not= 0$ and $l<k$ implies $\Coff [i,k]
b {\nu,l} = 0$ the sum \GSSumResT [] {j,k} b {\eta,l} \Coff \nu must
be nonzero. So there is $\eta \in \ProLevel j \nu$ with $\Coff [j,k] b
{\eta,l} \not= 0$. This contradicts the facts $l < k$ and \GSRes [j,k] b
k equals 0.
 \end{Proof}
\end{LEMMA}

\begin{LEMMA} {Res_eq}
 Suppose $\GS b \in \Gr$ and $I$ is a subset of \Set {(i,j) \in \Ind
[\gamma]} {\GSRes [i,j] b {\Braces j} = \GS [i,j] b}. Then for all
$(i,j) \in I$, $\nu \in \Level i$, and $k \in \Pro I i \Inter \Pro I
j$,
 \[
	\Coff [i,j] b {\nu,j}
	= \GSSumResT [] {j,k} b {\eta,k} \Coff \nu
	= \Coff [i,k] b {\nu,k}.
 \]

\begin{Proof}%
Since $(i,k)$ and $(j,k)$ are in $I$, both \Coff [i,k] b {\nu,j} and
\Coff [j,k] b {\eta,l} are equal to 0 for all $\eta \in \Level j$ when
$l \not= k$. Hence \EquationOfLemma {Eq} {=j} can be reduced to the
form $\Coff [i,j] b {\nu,j} = \GSSumResT [] {j,k} b {\eta,k} \Coff
\nu$. Now $(i,j) \in I$ guarantees that $\Coff [i,j] b {\nu,k} = 0$.
Thus the reduced form together with \EquationOfLemma {Eq} {>j} \Note
{applied for $l = k$} yield $\Coff [i,j] b {\nu,j} = \Coff [i,k] b
{\nu,k}$.
\end{Proof}
\end{LEMMA}

\begin{LEMMA} {Res_e.c.}
 Suppose \GS b is an element of \Gr.
 \begin{ITEMS}

\ITEM {Not_0}
 If \GS b is not in \Fact and $I$ is an eventually coherent subset of
\Ind such that $\GS [i,j] b = \GSRes [i,j] b {\Braces j}$ for all
$(i,j) \in I$, then there is an eventually coherent subset $J$ of $I$
with $\GS [i,j] b = \GSRes [i,j] b {\Braces j} \not= 0$ whenever
$(i,j) \in J$.

\ITEM {Bound} \Note{$\Cf \gamma > \aleph_0$}\ 
 If $J$ is an eventually coherent subset of \Ind such that $\GS [i,j]
b = \GSRes [i,j] b {\Braces j} \not= 0$ for all $(i,j) \in J$, then
there are a bound $n^* < \omega$ and an eventually coherent subset $K$
of $J$ such that $\Card {\Supp [i,j] b} < n^*$ for all $(i,j) \in K$.

\end{ITEMS}

\begin{Proof}
 \ProofOfItem {Not_0}
 Since \NotFactEquiv {\GS b} \Zero it follows by \Lemma
{Eventually_coherent} that there is no subset of $\Set {(i,j) \in I}
{\GS [i,j] b = 0}$ which would be eventually coherent. Hence there is
an unbounded subset $S$ of \First I such that for each $i \in S$ there
is $j_i \in \Pro I i$ with $\GS [i,j_i] b$ nonzero. Fix any $i \in S$.
Since $\GS [i,j_i] b = \GSRes [i,j_i] b {\Braces {j_i}} \not= 0$, let
$\nu_i$ be an element of \Level i with $\Coff [i,j_i] b {\nu_i,j_i}
\not= 0$. By \Lemma {Res_eq}, $\Coff [i,k] b {\nu_i,k} = \Coff
[i,j_i] b {\nu_i,j_i} \not= 0$ for all $k \in \Pro I i \Inter \Pro I
{j_i}$. Because $I$ was eventually coherent, we have shown that $J = I
\Inter (S \times S)$ is an eventually coherent set as wanted in the
claim.

 \ProofOfItem {Bound}
 First of all we claim that for each $i \in \First J$ the union
\BigUnion [j \in \Pro J i] {\Supp [i,j] b} is of finite cardinality.
Observe that for every $(i,j) \in J$, $\Supp [i,j] b = \Supp [i,j] b
\Inter (\Level i \times \Braces j)$.

Assume, contrary to this subclaim, that $i \in \First J$, \Seq {j_m}
{m < \omega} is an increasing sequence of ordinals in \Pro J i, and
\Set {\nu_m} {m < \omega} is a set of distinct elements from \Level i
such that \Coff [i,j_m] b {\nu_m,j_m} nonzero for every $m <
\omega$. Since $J$ is eventually coherent and $\gamma$ is of
uncountable cofinality let $k < \gamma$ be the minimal element in
$\Pro J i \Inter \BigInter [m < \omega] {\Pro J {j_m}}$. Now for each
$m < \omega$, the pairs $(i,j_m)$, $(i,k)$, and $(j_m,k)$ are in $J$,
and by \Lemma {Res_eq}, the equation $\Coff [i,j_m] b {\nu_m,j_m} =
\Coff [i,k] b {\nu_m,k} \not= 0$ holds. So the infinite set \Set
{(\nu_m,k)} {m < \omega} is a subset of \Supp [i,k] b, a
contradiction.

It follows from the subclaim that for each $i \in \First J$, the
finite ordinal
 \[
	n_i = \Card [\big] {%
	\BigUnion [j \in \Pro J i] {\Supp [i,j] b}} +1
 \]
 satisfies $\Card {\Supp [i,j] b} < n_i$ for all $j \in \Pro J
i$. Since \First J is uncountable, there are $n^* < \omega$ and an
unbounded subset $S$ of \First J such that $n_i = n^*$ for all $i \in
S$. So $n^*$ and the set $K = J \Inter (S \times S)$ meet the
requirements of the claim.
 \end{Proof}
\end{LEMMA}

\begin{LEMMA} {Step} \Note {$\Cf \gamma > \Card R$}\ 
 Suppose \GS b is in \Gr and $I$ is an eventually coherent subset of
\Set {(i,j) \in \Ind [\gamma]} {\GSRes [i,j] b {\Braces j} = \GS [i,j]
b \not= 0}. Then there are $d \in R$, $t \in \Br$, and an eventually
coherent subset $J$ of $I$ for which $\Coff [i,j] b {t(i),j} = d \not=
0$ whenever $(i,j) \in J$.

\begin{Proof}
 We define by induction on $\alpha < \Cf \gamma$ the following
objects:
 \begin{itemize}

\item an increasing sequence \Seq {i_\alpha} {\alpha < \Cf \gamma} of
ordinals in \First I with limit $\gamma$;

\item an increasing sequence $\Seq {\nu_\alpha} {\alpha < \Cf \gamma}
\in \Prod [\alpha < \Cf \gamma] {\Level {i_\alpha}}$;

\item subsets $K_\alpha$ of $\Pro I {i_\alpha}$ such that $\First I
\Minus K_\alpha$ are bounded in $\gamma$;

\item elements $d_\alpha \in R \Minus \Braces 0$ such that for every
$k \in K_\alpha$, $\Coff [i_\alpha,k] b {\nu_\alpha,k} = d_\alpha$.

\end {itemize}

This suffices since $\Card R < \Cf \gamma$ implies that there are $d
\in R$ and $H \Subset \Cf \gamma$ unbounded in \Cf \gamma such that
$d_\alpha = d$ for every $\alpha \in H$. Moreover, the claim is
satisfied by $t \in \Br$ and $J \Subset I$ defined as follows. For
every $i < \gamma$, $t(i) = \Res {\nu_{\beta_i}} i$, where $\beta_i =
\min \Set {\alpha < \Cf \gamma} {i_\alpha \geq i}$, and $J = \BigUnion
[\alpha \in H] {\Par {\Braces {i_\alpha} \times (S \Inter
K_\alpha)}}$, where $S$ is \Set {i_\alpha} {\alpha \in H}.

Let \Seq {\gamma_\alpha} {\alpha < \Cf \gamma} be an increasing
sequence with limit $\gamma$. Define
 \[
	i_\alpha = \Min [\big] {
	(\First I \Inter \BigInter [\beta < \alpha] {K_\beta})
	\Minus \gamma_\alpha}
 \]
 and
 \[
	j = \Min [\big] {
	\First I \Inter \Pro I {i_\alpha}
	\Inter \BigInter [\beta < \alpha] {\Pro I {i_\beta}}},
 \]
 where both \BigInter [\beta < \alpha] {K_\beta} and \BigInter [\beta
< \alpha] {\Pro I {i_\beta}} are equal to $\gamma$ when $\alpha =
0$. This pair $(i_\alpha,j)$ is well-defined since $I$ is eventually
coherent, $\alpha < \Cf \gamma$, and when $\alpha > 0$, $\First I
\Minus K_\beta$ is bounded for each $\beta < \alpha$ by the induction
hypothesis.

If $\alpha = 0$, then $(i_0, j) \in I$ guarantees that $\GSRes [i_0,
j] b {\Braces j} = \GS [i_0,j] b \not= 0$. Hence we can find $\nu_0
\in \Level {i_0}$ with $\Coff [i_0,j] b {\nu_0,j} \not= 0$.

When $\alpha > 0$ we define elements $\eta_\beta \in \Pro {\Level
{i_\alpha}} {\nu_\beta}$ for each $\beta < \alpha$ as follows. Fix
$\beta < \alpha$. Since $i_\alpha \in K_\beta$ we get by the induction
hypothesis that $\Coff [i_\beta,i_\alpha] b {\nu_\beta,i_\alpha} =
d_\beta \not= 0$. Furthermore $(i_\beta, i_\alpha) \in I$ \Note
{because $K_\beta \Subset \Pro I {i_\beta}$}, $(i_\beta,j) \in I$, and
$(i_\alpha,j) \in I$ together with \Lemma {Res_eq} yield
 \[
	\GSSumResT [] {i_\alpha,j} b {\eta,j} \Coff {\nu_\beta}
	= \Coff [i_\beta,i_\alpha] b {\nu_\beta,i_\alpha}	
	\not= 0.
 \]
 Therefore we can find $\eta_\beta \in \ProLevel {i_\alpha}
{\nu_\beta}$ for which $\Coff [i_\alpha,j] b {\eta_\beta,j} \not= 0$.

If $\alpha > 0$ is a successor ordinal define $\nu_\alpha$ to be
$\eta_{\alpha-1}$. When $\alpha$ is a limit ordinal, the finiteness of
the support \Supp [i_\alpha,j] b ensures that there are $\nu_\alpha
\in \Level {i_\alpha}$ and an unbounded subset $H$ of $\alpha$ such
that $\eta_{\beta'} = \nu_\alpha$ for all $\beta' \in H$. By the
induction hypothesis $\nu_\beta \Tleq \nu_{\beta'}$ for all $\beta <
\beta' < \alpha$. Hence $\nu_\beta \Tleq \nu_{\beta'} \Tleq
\eta_{\beta'} = \nu_\alpha$ holds for every $\beta < \alpha$ and
$\beta' = \Min {H \Minus \beta}$.

Let $d_\alpha$ be \Coff [i_\alpha,j] b {\nu_\alpha,j}. By \Lemma
{Res_eq}, every $k \in \Pro I {i_\alpha} \Inter \Pro I j$ satisfies
that $\Coff [i_\alpha, k] b {\nu_\alpha, k} = \Coff [i_\alpha, j] b
{\nu_\alpha, j} = d_\alpha$. Hence $i_\alpha$, $\nu_\alpha$, and
$d_\alpha$ together with the set $K_\alpha = \Pro I {i_\alpha} \Inter
\Pro I j$ meet the requirements given at the beginning of the proof.
 \end{Proof}
\end{LEMMA}

\begin{COROLLARY} {Empty_Br} \Note{$\Cf \gamma > \aleph_0$}\ 
 If $\Br$ is empty, then $\Gr = \Fact$.

\begin{Proof}
 Suppose $\GS a \in \Gr \Minus \Fact$. By \ItemOfLemma {Res}
{Cobounded} together with \ItemOfLemma {Res_e.c.} {Not_0} there is
$\GS b \in \Gr$ such that \FactEquivGS a b and the set \Set {(i,j) \in
\Ind [\gamma]} {\GSRes [i,j] b {\Braces j} = \GS [i,j] b \not= 0} is
eventually coherent. By \Lemma {Step} there is a $\gamma$-branch
through the tree $T$, i.e., $\Br \not= \emptyset$. Observe that the
assumption $\Card R < \Cf \gamma$ is not needed, as can be seen from
the proof of \Lemma {Step}.
 \end{Proof}
\end{COROLLARY}

\begin{LEMMA} {Generate}%
 \Note{$\Cf \gamma > \max \Braces {\aleph_0, \Card R}$}\ 
 The elements $\GenS t$, $t \in \Br$, generate \Gr modulo \Fact.

\begin{Proof}
 We show that for every $\GS a \in \Gr$ with $\GS a \not\in \Fact$ we
can find $n < \omega$, $d_1, \dots, d_n \in R \Minus \Braces 0$
and $t_1, \dots, t_n \in \Br$ satisfying
 \EQUATION{Goal}{
	\FactEquiv {\GS a} 
	{\Sum [1 \leq m \leq n] {d_m \Times \GenS {t_m}}}.
 }%EndEquation

Suppose $\GS a \in \Gr \Minus \Fact$. By \ItemOfLemma {Res}
{Cobounded} and \ItemOfLemma {Res_e.c.} {Not_0} let \GS b be an
element of \Gr and $I_1$ an eventually coherent subset of \Ind such
that \FactEquivGS a b and for each $(i,j) \in I_1$, $\GS [i,j] b =
\GSRes [i,j] b {\Braces j} \not= 0$. Furthermore, we may assume by
\ItemOfLemma {Res_e.c.} {Bound} that $n^* < \omega$ is a bound for
which $\Card {\Supp [i,j] b} < n^*$ hold for all $(i,j) \in I_1$.

By \Lemma {Step} there are $d_1 \in R$, $t_1 \in \Br$, and an
eventually coherent set $J_1 \Subset I_1$ having the property that
$\Coff [i,j] b {t_1(i), j} = d_1 \not= 0$ whenever $(i,j) \in
J_1$. Since $d_1 \Times \GenS {t_1} \in \Gr$, the sequence $\GS c =
\GS b - d_1 \Times \GenS {t_1}$ is in \Gr. If \GS c is in \Fact, then
\FactEquiv {\GS b} {d_1 \Times \GenS {t_1}}, and because of
\FactEquivGS a b, also \Equation {Goal} holds for $n=1$.

Suppose $1 \leq n < \omega$ and objects $d_m \in R \Minus \Braces 0$,
$t_m \in \Br$, and $J_m \Subset J_1$ for $m \leq n$ are already
defined. Assume also that these objects satisfy the following
conditions:
 \begin{enumerate}

\ITEM {chain}
 $J_{m'} \Superset J_m$ for all $1 \leq m' \leq m \leq n$;

\ITEM {diff_t}
 for all $1 \leq m' < m \leq n$ and $i \in \First {(J_m)}$,
 $t_{m'}(i) \not= t_m(i)$;

\ITEM {d_m}
 for every $1 \leq m \leq n$ and $(i,j) \in J_m$, $\Coff [i,j] b
{t_m(i), j} = d_m \not= 0$;

\ITEM {c}
 $\GS c = \GS b - \Sum [1\leq m \leq n] {d_m \Times \GenS {t_m}}
\not\in \Fact$.

\end{enumerate}

Clearly $\GS [i,j] c = \GSRes [i,j] c {\Braces j}$ and $\Card {\Supp
[i,j] c} \leq \Card {\Supp [i,j] b} < n^*$ for all $(i,j) \in
J_n$. Again by \ItemOfLemma {Res_e.c.} {Not_0}, there is an eventually
coherent set $I_{n+1} \Subset J_n$ such that for each $(i,j) \in
I_{n+1}$, $\GS [i,j] c \not= 0$. Moreover, by \Lemma {Step}, there are
$d_{n+1} \in R$, $t_{n+1} \in \Br$, and an eventually coherent set
$J_{n+1} \Subset \Set {(i,j) \in I_{n+1}} {\Coff [i,j] c {t_{n+1}(i),
j} = d_{n+1} \not= 0}$.

The properties \Item {diff_t}, \Item {d_m} and \Item c above imply
that $\Coff [i,j] c {t_m(i),j} = \Coff [i,j] b {t_m(i),j} - d_m = 0$
for every $m \leq n$ and $(i,j) \in J_m$. On the other hand, \Coff
[i,j] c {t_{n+1}(i), j} is nonzero for each $(i,j) \in J_{n+1}$. Thus
$t_{n+1}(i)$ can not be in \Set {t_m(i)} {1\leq m \leq n} if $i \in
\First {(J_{n+1})}$. So for all $(i,j) \in J_{n+1}$, $\GenT [i,j]
{t_{n+1}} \not\in \Set {\GenT [i,j] {t_m}} {1\leq m \leq n}$, and
consequently $\Coff [i,j] b {t_{n+1}(i), j} = \Coff [i,j] c
{t_{n+1}(i), j}$. Thus also $J_{n+1}$, $t_{n+1}$, and $d_{n+1}$
satisfy the properties \Item {chain}, \Item {diff_t}, and \Item {d_m}
\Note {but not necessarily \Item c}.

We claim that there must be $n < n^*$ such that
 \EQUATION{Subclaim}{
	\GS b -
	\Sum [1 \leq m \leq n] {d_m \Times \GenS {t_m}} \in \Fact.
 }%EndEq

Assume, contrary to this subclaim, that the process introduced above
has been carried out $n^*$ many times and objects $J_m$, $t_m$, $d_m$
for $i \leq m \leq n^*$ are defined. In addition to that suppose they
satisfy the conditions \Item {chain}, \Item {diff_t}, and \Item
{d_m}. Define $i = \Min [\big] {\First{(J_{n^*})}}$ and $j = \Min
{\Pro {J_{n^*}} i}$. Then for every $m \leq n^*$, $(i,j) \in J_m$
yields $\Coff [i,j] b {t_m(i),j} = d_m \not= 0$. This contradicts the
condition $\Card {\Supp [i,j] b} < n^*$, since the set $\Set
{(t_m(i),j)} {m \leq n^*} \Subset \Supp [i,j] b$ is of cardinality
$n^*$.

Now suppose $n < \omega$ is a finite ordinal satisfying \Equation
{Subclaim}. Then \FactEquiv {\GS b} {\Sum [1 \leq m \leq n] {d_m
\Times \GenS {t_m}}}, and because \FactEquivGS a b also \Equation
{Goal} is satisfied.
 \end{Proof}
\end{LEMMA}

\begin{CONCLUSION}{System}
 For any ordinal $\gamma$ of uncountable cofinality, ring $R$ with
$\Card R < \Cf \gamma$, and tree $T$ of height $\gamma$, the inverse
$\gamma$-system $\System = \SystemSeq$ has the properties that
 \[
	\Card {G_i} = \max \Braces {%
		\Card \gamma, \Card {\Level i}, \Card R
	}
 \]
 for all $i < \gamma$, and
 \[
	\CardGrF = \Card [\big] \BrGen.
 \]
 \end{CONCLUSION}

\begin{RefArea}{Env}{Theorem}{System}

\begin{SeparateProof}
 Remember that $\lambda$ and $\kappa$ were cardinals with $\aleph_0 <
\kappa = \Cf \lambda < \lambda$. We wanted to study possible
cardinalities $\mu$ of the quotient limit \GrF [A], where \System [A]
is an inverse $\kappa$-system consisting of abelian groups having
cardinality $< \lambda$. Now \Conclusion {System} gives a complete
solution to this problem because of $\lambda > \Cf \lambda = \kappa =
\Cf \kappa > \aleph_0$. Namely, in order to meet the requirements
$\Card {G_i} < \lambda$ for all $i < \kappa$, it is needed only to
ensure that $R$ and the \Th i level of $T$ are small enough. On the
other hand, a suitable choice of $R$ and $T$ yields any desired value
for $\mu = \CardGrF$. We briefly describe methods to choose suitable
$R$ and $T$ for every nonzero $\mu \leq \lambda^\kappa$.

For any $R$, \CardGrF equals 1 when \Br [\kappa] is empty. So $\mu=1$
is possible since obviously there exists a tree of height $\kappa$
without $\kappa$-branches and having levels of cardinality $< \lambda$
when $\lambda$ singular of cofinality $\kappa$. Also all the finite
values $\mu > 1$ are possible by taking $T$ with only one
$\kappa$-branch and $R$ with $\Card R = \mu$.

Furthermore the case of infinite $\mu < \lambda$ is satisfied by any
$R$ with $\Card R < \min \Braces {\kappa, \mu}$ and $T$ with exactly
$\mu$ many $\kappa$-branches. The value $\mu = \lambda$ is possible
for any $R$ with $\Card R < \kappa$ because a suitable tree can be
constructed, for example, as follows. Let \Seq {\lambda_i} {i <
\kappa} be an increasing sequence of ordinals $< \lambda$ with limit
$\lambda$. Then the tree
 \[
	T = \Set {\Res t \alpha}
	{\alpha < \kappa, t \in \Prod [i < \kappa] {\lambda_i}, \And
	t(i) \Text{is nonzero only for finitely many} i < \kappa},
 \]
 ordered by inclusion, satisfies $\Card {\Br [\kappa]} = \lambda$ and
$\Card {\Level i} = \lambda_i < \lambda$ for each $i < \kappa$.

Also the cardinalities $\mu$ of the quotient limit, when $\lambda <
\mu \leq \lambda^\kappa$, are possible for any ring of cardinality $<
\kappa$. Existence of a suitable tree is proved for example in \cite
[Fact 10] {Sh262} under the assumption that $2^\kappa < \lambda$ and
$\theta^{<\kappa} < \lambda$ for every $\theta < \lambda$ \Note {other
sources for a proof are given in \cite [Analytical Guide \S 10]
{Sh:g}}.
 \end{SeparateProof}

\end{RefArea}

%end(Section_System)=========================================

\end{SECTION}

%\InputText{}{BibAndAddresses}
%begin(BibAndAddresses)=======================================

%\bibliographystyle{alpha}
%\bibliography{abbreviations,citations,mathscinet}

%begin(SHVa:644.bbl)=======================================

%end(SHVa:644.bbl)=========================================
%-----------------------------------------------------------------------
\catcode`@=11
\begin{tabbing}
Saharon Shelah:\= \\
	\>Institute of Mathematics\\
	\>The Hebrew University\\
	\>Jerusalem. Israel\\
\\
	\>Rutgers University\\
	\>Hill Ctr-Busch\\
	\>New Brunswick. New Jersey 08903\\
	\>\texttt{shelah@math.rutgers.edu}
\end{tabbing}

\begin{tabbing}
Pauli V\"{a}is\"{a}nen:\= \\
	\>Department of Mathematics\\
	\>P.O. Box 4\\
	\>00014 University of Helsinki\\
	\>Finland\\
	\>\texttt{pauli.vaisanen@helsinki.fi}
\end{tabbing}
\catcode`@=13

%-----------------------------------------------------------------------

%end(BibAndAddresses)=========================================

\end{document}